\documentclass[11pt]{article} 
	 
\usepackage{amsmath} 
\usepackage{psfrag,graphicx,verbatim,array,multicol,palatino,enumerate}
\usepackage{amsmath,amsfonts,bm,rotating,natbib}
\usepackage{pstricks,pst-node}
\usepackage{epsfig,setspace}

\def\bz{{\rm {\bf{z}}}}
\def\bw{{\rm {\bf{w}}}}
\def\ba{{\rm {\bf{a}}}}

\def\bB{{\rm {\bf{B}}}}

\def\bN{{\rm {\bf{N}}}}
\def\bY{{\rm {\bf{Y}}}}
\def\btheta{\boldsymbol{\theta}}
\def\bTheta{\boldsymbol{\bTheta}}

\def\bTheta{\boldsymbol{\Theta}}
\def\bdelta{\boldsymbol{\delta}}
 
\begin{document}

\begin{center} 
 
{\large{\bf BAYESIAN CLUSTERING OF TRANSCRIPTION FACTOR BINDING MOTIFS}}

\bigskip

SHANE T. JENSEN\footnote{Department of Statistics, The Wharton School, University of Pennsylvania, Philadelphia, PA 19104  \\ {\tt{stjensen@wharton.upenn.edu}}} and 
JUN S. LIU\footnote{Department of Statistics, Harvard University, Cambridge, MA 02138 \\  {\tt{jliu@stat.harvard.edu}}}. 
 
\end{center}

\bigskip

\begin{abstract} 
 
Genes are often regulated in living cells by proteins called
transcription factors (TFs) that bind directly to short segments of
DNA in close proximity to specific genes.  These binding sites
have a conserved nucleotide appearance, which is called a motif.  Several recent studies of transcriptional regulation require the reduction of a large collection of motifs into clusters based on the similarity of their nucleotide composition.       
We present a principled approach to this clustering problem based upon a Bayesian hierarchical model that accounts for both within- and 
between-motif variability. We use a Dirichlet process 
prior distribution that allows the number of clusters to vary and we also present a novel generalization that allows the core width of each motif to vary.     

This clustering model is implemented, using a Gibbs sampling strategy,
on several collections of transcription factor motif matrices.   Our stochastic
implementation allows us to examine the variability of our results in addition to focusing on a set of best clusters. 
Our clustering results identify several motif clusters that suggest several transcription factor protein families are actually mixtures of several smaller groups of highly similar motifs, which provides substantially more refined information compared with the full set of motifs in the family.  
Our clusters provide a means by which to organize transcription factors based on binding motif similarities, which can be used to reduce motif redundancy within large databases such as JASPAR and TRANSFAC, which aides the use of these databases for further motif discovery. 
Finally, our clustering procedure has been used in combination with discovery of evolutionarily-conserved motifs to predict co-regulated genes.  
An alternative to our Dirichlet process prior distribution is presented that differs substantially in terms of a priori clustering characteristics, but shows no substantive difference in the clustering results for our dataset.  
Despite our specific application to transcription factor binding motifs, our Bayesian clustering model based on the Dirichlet process has several advantages over traditional clustering methods that could make our procedure appropriate and useful for many clustering applications.  Software for our method is available at {\tt http://stat.wharton.upenn.edu/$\sim$stjensen/research/cluster.html} .   

\end{abstract} 

\bigskip

\section{Introduction}\label{introduction}

The complete information that defines the characteristics of living
cells within an organism is encoded in the form of a moderately
simple molecule, deoxyribonucleic acid,  or DNA.  The building 
blocks of DNA are four nucleotides, abbreviated by their attached
organic bases as A, C, G, and T.  A-T and C-G
are complementary bases between which hydrogen bonds can form.  A DNA
molecule consists of two long chains of nucleotides that are
complimentary to each other and joined by hydrogen bonds twisted into
a double helix. This structure gives rise to the term {\it base pair}
when describing a DNA sequence.  A long DNA chain in a living
cell is called a {\it chromosome}.  For example, every human cell has 23
pairs of chromosomes with lengths ranging from 47 million  to  245
million base pairs (Mb). The specific ordering of the four types of
nucleotides along these chains is the means by which the information
is  stored that completely defines all functions within a cell. The
term {\it gene} refers to sequence segments along a chromosome that are
used to code the information for making {\it proteins}, the fundamental action
molecules of the cell. Surprisingly, only about 1 to 2
\% of the entire human genome (set of chromosomes) corresponds to gene regions. It is believed
that much of the mechanism for controlling when, where, and how
much protein will be produced is located in the ``non-coding''
region located {\it upstream} (ie. directly before) the gene sequence.  

Transcribing or activating a gene requires not only 
the DNA sequence in the upstream region, but also many proteins called
transcription factors (TF). When these TFs are present, 
they bind to specific DNA patterns in the upstream sequence of genes, and either 
induce or repress the 
transcription of these genes by recruiting other necessary 
proteins \cite[]{LodBalBerZipMatDar1995}.  A particular transcription factor protein is able to bind and regulate only certain target genes by recognizing a short (6-20 basepairs long) sequence of nucleotides called a transcription factor binding site (or, more simply, a {\it site}).  Different binding sites (located near different genes) of the same transcription factor protein show a substantial sequence conservation, which we call a {\it motif}, but some variability is also present.  

\subsection{Statistical Formulation of Motifs} \label{motifformulation}

Each motif is mathematically formulated as a {\it motif matrix},  
which measures the desirability of each base at each
position of the motif.  The simplest matrix is an alignment or count matrix
$\bN_{jk}$, which records the occurrence of  base $k$ at position $j$ of 
all the sites for this motif.   Table~\ref{matrix} shows the count matrix for motif MA0011 from the database we will use in Section~\ref{jasparapp}.  
Also shown
in Table~\ref{matrix} is  
the corresponding frequency matrix ($f_{jk} = \bN_{jk}/N$, where $N$ is
the number of  motif sites) for motif MA0011.   

\begin{table}[ht]
\caption{Matrix representations of the motif MA0011}\label{matrix}
\begin{center}
{\footnotesize
\begin{tabular}{ccccccccccccccccccc}
\multicolumn{9}{c}{\bf{Count Matrix}} &  & \multicolumn{9}{c}{\bf{Frequency Matrix}} \\
Pos & 1 & 2 & 3 & 4 & 5 & 6 & 7 & 8 & \hspace{0.5cm} & Pos & 1 & 2 & 3 & 4 & 5 & 6 & 7 & 8 \\
A  &  3  & 5 & 0 & 0  & 12  & 1 & 2 & 1 & \hspace{0.5cm} & A  &  0.25  & 0.42 & 0.00 & 0.00  & 1.00  & 0.08 & 0.17 & 0.08\\
C  & 1 & 2  & 10  & 1 & 0 & 1 & 0  & 2 & \hspace{0.5cm} & C  & 0.08 & 0.17  & 0.83  & 0.08 & 0.00 & 0.08 & 0.00  & 0.17 \\
G  & 1 & 1 & 0 & 0 & 0 & 2 & 1 & 1 & \hspace{0.5cm} & G  & 0.08 & 0.08 & 0.00 & 0.00 & 0.00 & 0.17 & 0.08 & 0.08 \\
T  & 7 & 4 & 2  & 11 & 0 & 8 & 9  & 8 & \hspace{0.5cm} & T  & 0.58 & 0.33 & 0.17  & 0.92 & 0.00 & 0.67 & 0.75  & 0.67 \\
\end{tabular}
}
\end{center}
\end{table}

\begin{figure}[ht]
\caption{Sequence logo of the motif MA0011}\label{seqlogo}
\begin{center}
\hspace{-0.5cm}
\includegraphics[width=4in,height=1in]{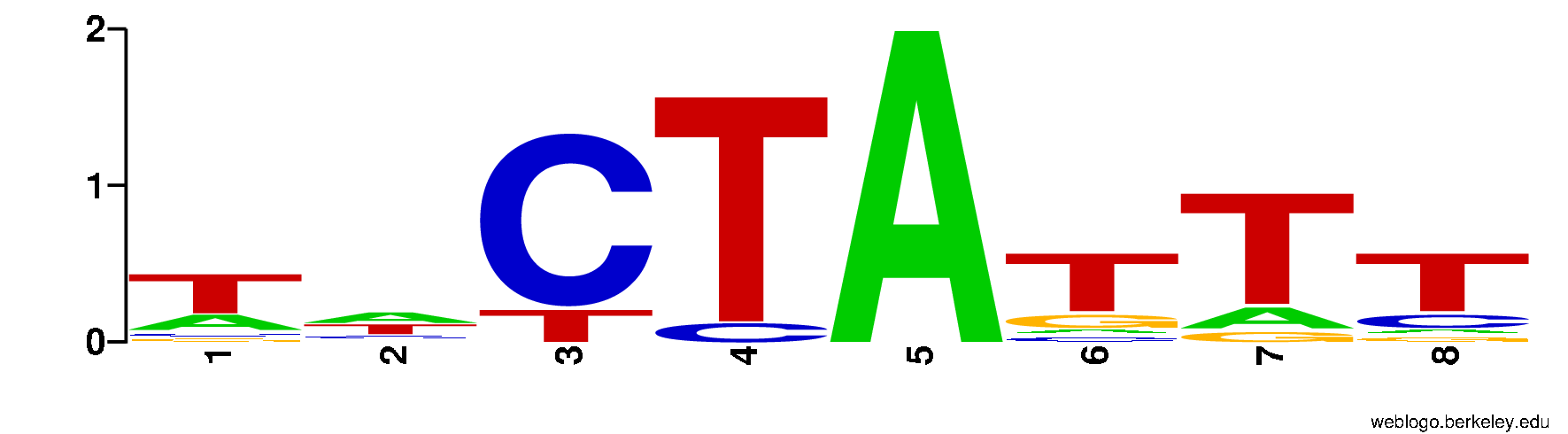}
\end{center}
\end{figure}

\cite{SchSte1990} used the motif matrix to construct a {\it Sequence Logo} as a 
means by which to visualize the appearance of the motif.  Figure~\ref{seqlogo} gives the sequence 
logo for the same motif {\tt MA0011}, as constructed by the program {\it WebLogo} \cite[]{weblogo}.  The height of each 
position is equal to its {\it information content} ($\sum_k f_{jk} \log [ f_{jk}/\theta_{0k} ]$  where $\theta_{0k}$ is the 
proportion of base $k$ in the non-motif background positions) and the size of each letter 
is proportional to the letter's relative frequency $f_{jk}$.

\subsection{Discovery of Motifs}\label{motifdiscoveryintro}

As reviewed in \cite{JenLiuZhoLiu2004}, Bayesian motif discovery models are the foundation of many popular programs for discovering conserved binding sites in large sequence datasets.   These models are usually implemented by an iterative strategy, such as the EM algorithm \cite[]{DemLaiRub1977} or the Gibbs sampler \cite[]{GemGem1984}, that utilizes the ease with which the motif frequency matrix can be estimated if the binding site locations are known, and the corresponding ease with which the binding site locations can be estimated if the motif appearance is known.    \cite{JenLiuZhoLiu2004} also discuss several constraints of these motif discovery procedures, such as assumptions of known motif width (the number of columns in the motif frequency matrix) and known abundance of binding sites.
\cite{JenLiu2004} demonstrates that allowing the motif width to vary leads to more accurate results in several motif discovery applications.  
 
Although the discovery and characterization of a single motif is often the 
goal of a particular biological investigation (see, for example, \cite{sigmaE}), it is common for scientists to be interested in examining the similarities and differences between an entire collection of discovered TF motifs.   Large collections of discovered motifs have been utilized in various applications, such as the phylogenetic discovery of co-regulated genes \cite[]{QinMccThoMayLawLiu2003} and the prediction of synergistic relationships between transcription factors 
\cite[]{HanLev2002}.   These applications each represent a highly specialized approach to the utilization of a collection of motifs and do not address the issue of a 
general statistical approach for the sharing of information between motifs.  

\subsection{Modeling Motif Similarity by Clustering}\label{clusintro}

\begin{figure}[ht]
\caption{Four different motifs from the JASPAR database} \label{differentmotifs}
\begin{center}
\begin{tabular}{ccc}
MA0011 & \hspace{1cm} & MA0015 \\
\includegraphics[width=2.75in,height=1in]{fig.logo.MA0011.ps} 
 & & 
\includegraphics[width=2.75in,height=1in]{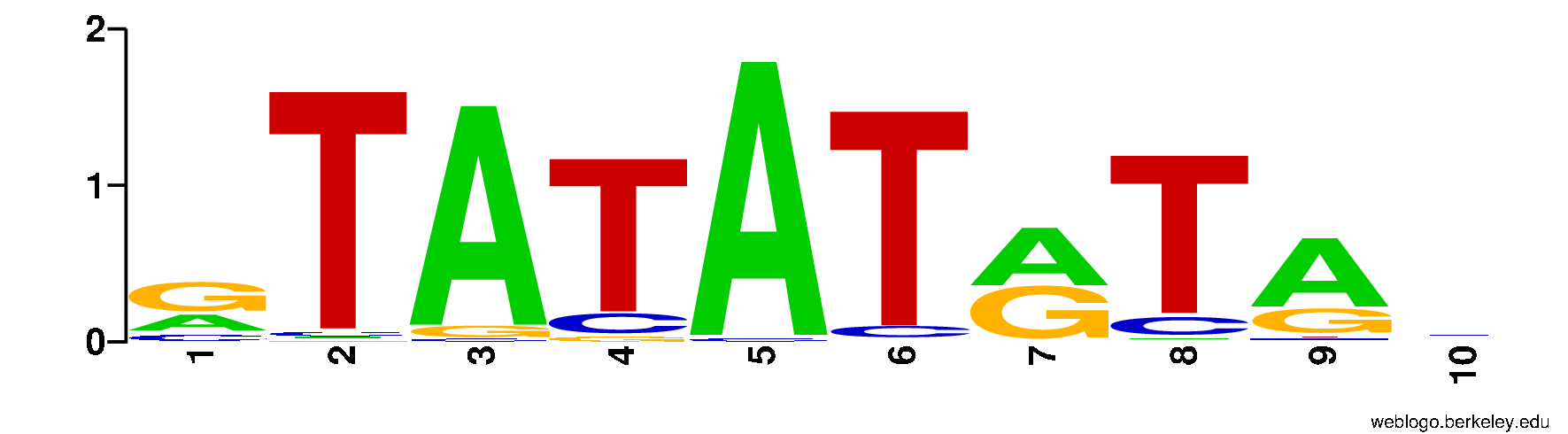} 
 \\
\phantom{nothing} & & \\
MA0031 & & MA0032\\
\includegraphics[width=2.75in,height=1in]{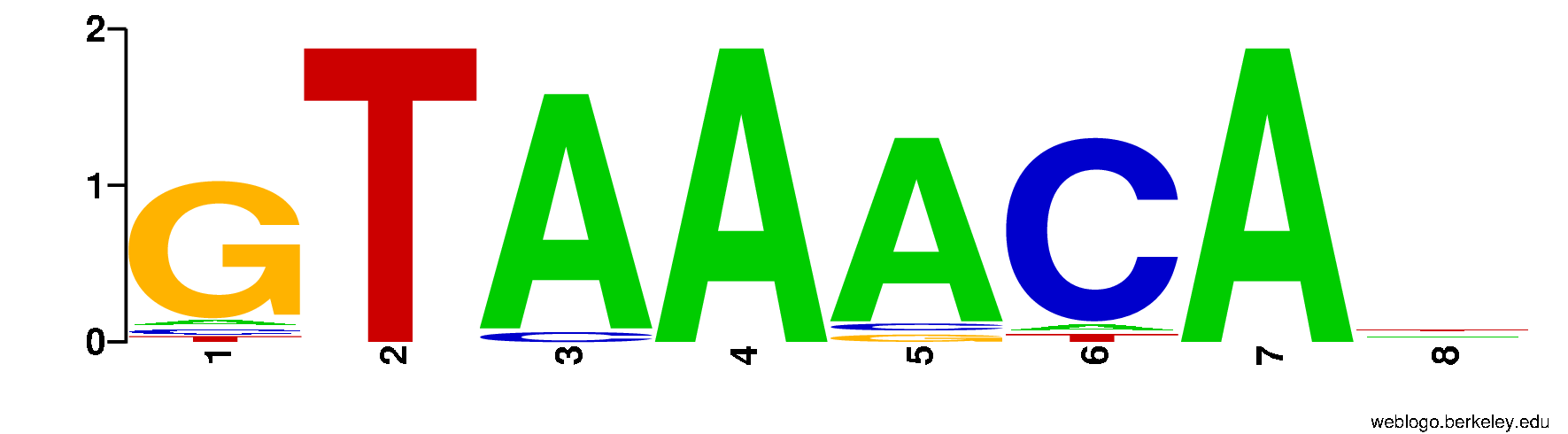} 
 & & 
\includegraphics[width=2.75in,height=1in]{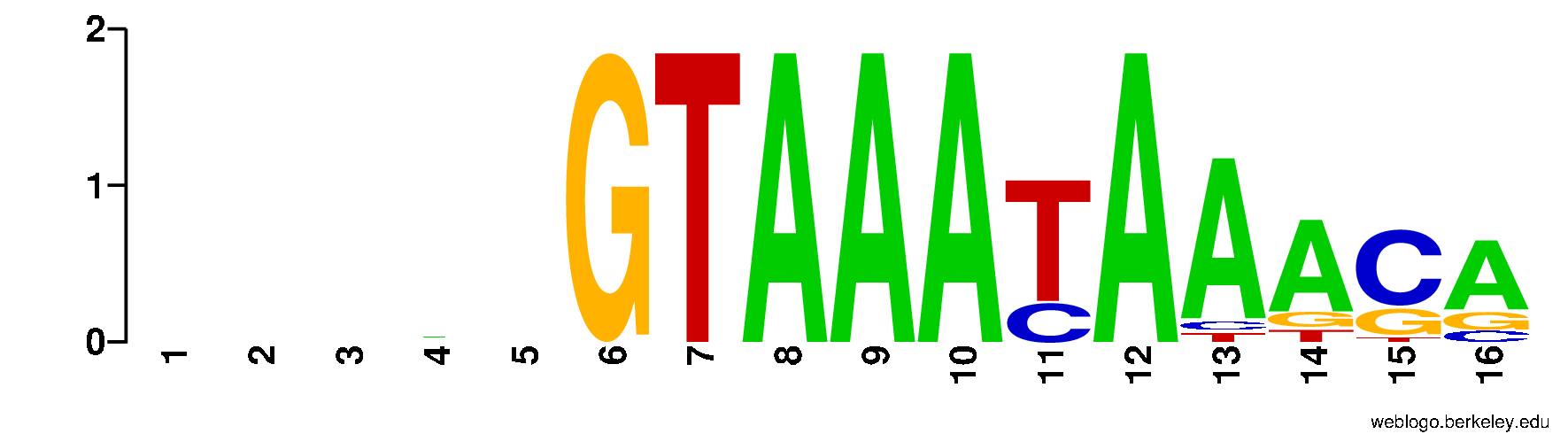} 
 \\
\end{tabular}
\end{center}
\end{figure}

Figure~\ref{differentmotifs} 
shows the sequence logos for four different motifs from the JASPAR database which will be analyzed in Section~\ref{jasparapp}.
It is clear that these motifs all show differences in width and appearance, but there exists some similarity or 
common structure within this set. One could certainly 
argue for grouping MA0031 and MA0032 together based upon their similar appearances, though they do differ in both width and composition , but this 
decision is based on ad-hoc personal judgement.
The statistical problem of interest here is to model the common structure
between these different motifs and find a principled means by which to group motifs together based upon their similarity.

\vspace{0.5cm}

There are several traditional statistical techniques for clustering observations 
together which are reviewed in \cite{Har1975}.  
{\it Hierarchical Tree Clustering} joins observations 
together into successively larger clusters based upon some sort of similarity measure.  
{\it K-means Clustering} groups observations into a pre-determined 
number of clusters by minimizing a within-cluster distance measure.  
Each of these techniques have elements that are not ideally suited for
our desired goal of motif clustering.  Hierarchical tree clustering
requires the user to specify a distance measure between the
observations (in this case, motif matrices), and it is not clear for
comparing motifs what type of simple distance metric should be used.  
In addition, the result of this algorithm is a tree that joins all
observations together, and it is not clear where the tree should be
``cut'' in order to produce a set of clusters.  \cite{Kielbasa2005} proposed two different distance metrics for the clustering of motif matrices, but needed to impose arbitrary thresholds on those distances in order to produce a set of clusters.   Similarly, \cite{SchSumZha2005} present several distance metrics for comparing motif matrices, and group together pairs of matrices below a p-value-based threshold.   K-means clustering is ideally suited for situations where the number of clusters is known {\it a
priori}.  When the number of clusters is unknown, K-means clustering becomes more difficult, and usually a 
cross-validation strategy is employed to estimate the number of clusters.  For our motif clustering applications, 
there is very little prior idea of how many motifs might cluster together in a particular collection
of motifs and so we seek a model that easily allows for an unknown number of clusters.  

In addition, these techniques consider the observations themselves as fixed and known, 
which is not the case for our applications where each
motif is only an estimate generated by a prior motif discovery procedure.
Recognizing that our discovered motifs themselves are estimated
quantities, we need to model both within-motif and between-motif
variability.  In Section~\ref{clusteringmodel}, we outline a Bayesian
hierarchical clustering model that encompasses both levels of
variability and does not require prior knowledge of the number of clusters.                                                                                                                                                              
As discussed briefly in Section~\ref{motifdiscoveryintro}, most motif discovery procedures assume a fixed and known motif width, but in reality 
the width is often unknown for many motifs and can vary substantially between different motifs.  In Section~\ref{varyingwidth}, we extend our Bayesian motif 
clustering model to allow each motif width to be an unknown variable that will be estimated by our procedure.  
Our model is implemented stochastically by a Gibbs sampling algorithm, which allows us to
examine not only ``best estimates" of motif clusters, but also the variability within our clustering results.
In Section~\ref{jasparapp}, 
we present various techniques for summarizing and understanding the results from our 
clustering procedure, within the context of an application to the JASPAR database of transcription factor binding motifs.  

These clusters provide a way to organize transcription factor proteins based on the similarity of their binding motifs.  We utilize this clustering to  analyze large databases of transcription factor matrices,  TRANSFAC \cite[]{transfac} and JASPAR \cite[]{jaspar} for motifs with high similarity (Section~\ref{combined}) and discuss how this clustering information can be used to refine the search for additional transcription factor binding sites.    In addition, we briefly discuss combining our clustering procedure with motif discovery in sequences that have been conserved by evolution, in order to predict genes that share similar transcription factor binding motifs and thus are possibly co-regulated (Section~\ref{crossclus}).

\section{Bayesian Motif Clustering Model}\label{clusteringmodel}

We use a Bayesian hierarchical model to infer common structure, in the
form of clusters, within a collection of discovered motifs. 
The data for each discovered 
motif is a count matrix $\bN_i$ which can have different widths and
number of counts (ie. number of binding sites included in the matrix) compared to other TF motifs.  
For now, we assume that our clustering will be
based on a motif matrices with a fixed and known width $w$, so we assume each of
these $n$ raw motif matrices  $\bN_i$ should contain a submatrix $\bY_i$ of dimension $w \times 4$ that will be considered the
{\it core} upon which the clustering will be based.  We will later extend our model to allow the core width within each cluster of motifs to vary.

\subsection{Hierarchical Framework}\label{hierarchicalmodel}

Hierarchical models are useful in a variety of scientific problems 
when the structure of the data suggests multiple levels of uncertainty.  
We want to include components for both within-motif and  
between-motif variability of the nucleotide counts  
$Y_{ijk}$ where $i$ indexes the motif, $j$ indexes the $w$ columns within  
each motif core, and $k$ indexes the four possible nucleotides ($k = a,c,g,$ or $t$) within each column.   
Our model on the within-motif variability between different binding sites for a count
motif $\bY_i$ is a product-multinomial model.  
We assume that each position (column) of the core count matrix $\bY_i$ follows an independent multinomial 
distribution parameterized by the corresponding column of an unknown frequency matrix $\bTheta_i$, ie.,

Within-motif level:  Count matrices $p(\bY_i | \bTheta_i) = \prod\limits_{j=1}^w p(\bY_{ij} | \btheta_{ij})$  
\begin{eqnarray*} 
\bY_{ij} = (Y_{ija},Y_{ijc},Y_{ijg},Y_{ijt}) \sim {\rm Multinomial} (n_i,\btheta_{ij}=(\theta_{ija},\theta_{ijc},\theta_{ijg},\theta_{ijt})) 
\end{eqnarray*} 

For our between-motif variability, we simply assume that each motif
frequency matrix $\bTheta_i$ in our collection share a common but
completely unknown distribution, denoted ${\rm F} (\cdot)$, ie.

Between-motif level: Frequency matrices $p(\bTheta_i)$ 
\begin{eqnarray*} 
\bTheta_i = (\btheta_{i1}, \ldots, \btheta_{iw}) \sim {\rm F} (\cdot)  
\end{eqnarray*} 
where ${\rm F} (\cdot)$ is an unknown distribution with $w$ dimensions
for the columns $\times$ 4 dimensions for the nucleotides (constrained
to sum to one).
This unknown distribution ${\rm F}(\cdot)$ represents the common
structure between the different motifs in the dataset.  Estimation of
this unknown distribution is complicated by the fact that our
frequency matrices $\bTheta_i$ are unknown, with only the count 
matrices $\bY_i$ being observed.
A popular Bayesian approach to non-parametric problems is to give the unknown distribution ${\rm F}(\cdot)$ a Dirichlet process prior, ${\cal D}(\gamma)$, where $\gamma$ is a finite (non-negative) measure, typically smooth \cite[]{Fer1974}. Here, since we have a multidimensional ${\rm F}(\cdot)$, we use a Dirichlet process prior  
${\cal D}(\gamma_1 \times \cdots \times \gamma_w)$, where each smooth measure $\gamma_j$ is 
four dimensional, taking the form of $\gamma_j = b \times {\rm Dirichlet}(\alpha,\ldots,\alpha)$ for  $j=1,\ldots,w$.  The parameter $b$ is a weighting factor that characterizes how close the unknown distribution $F$ is to the shape of $\gamma$ and the smoothness of $F$.  

\subsection{Clustering of Observations}

An important consequence of our model is that it enables 
similar motifs to be clustered together into a group modeled by one common frequency matrix.  
As explained in \cite{Fer1974}, if $\bTheta_1,\ldots,\bTheta_n$ are $n$ i.i.d. observations from the probability function $F$ whose prior distribution is the Dirichlet process ${\cal D}(\gamma)$, where $\gamma$ is a finite measure on the domain, then 
\begin{eqnarray*}  
{\rm F}(\cdot) | \bTheta_1,\ldots,\bTheta_n \sim {\cal D}(\gamma^\star) =  
 {\cal D}(\gamma + \sum\limits_{j=1}^n  \bdelta_{\bTheta_j} ) 
\end{eqnarray*} 
Thus, the posterior mean of ${\rm F}(\cdot)$, or the predictive distribution of a new observation, is proportional to $\gamma + \sum\limits_{j=1}^n  \bdelta_{\bTheta_j}$.  If the $\bTheta$'s only take on $C$ distinctive values, then we have a mixture of the smooth measure   
$\gamma$ and $C$ point masses (with potentially different weights).  These point mass components allow for the clustering of similar
observations. 
If we were to draw an additional $(n+1)$-th observation
$\bTheta^\star$ from this distribution ${\rm D}(\gamma^\star)$, that new observation
would either come from the smooth measure $\gamma$, or would take on a
value exactly equal to one of the current $\bTheta_j$'s, say $\bTheta_c$,
in which case $\bTheta_c$ and $\bTheta^\star$ are defined as being in the same
cluster.
The conditional distribution $p(\bTheta_i | \bTheta_{-i})$ of one current
observation $\bTheta_i$, given all other observations $\bTheta_{-i}$, is
also a mixture between the smooth measure and $C$ point masses at each
of the $\tilde{\bTheta}_{-i}$ that represent the unique values within
$\bTheta_{-i}$.  Any observations $\bTheta_m$ and $\bTheta_n$  that have the same value are defined as being in the
same cluster. This conditional distribution allows us to implement our model via 
Gibbs sampling \cite[]{GemGem1984}.   The Dirichlet process has been used as a prior distribution in nonparametric Bayesian analyses, such as \cite{Maceachern1994} and \cite{Esc1994} for the estimation of normal means and \cite{Liu1996} in a binomial hierarchical setting.  \cite{GreRic2001} 
discuss the use of the Dirichlet process as a flexible model for clustering observations, and 
present an extended class of Dirichlet-Multinomial allocations for which the Dirichlet process 
is a limiting case.  \cite{MedSiv2002} uses the clustering properties of the Dirichlet process 
prior as part of a hierarchical model for gene expression profiles from microarray data.   

\subsection{Gibbs Sampling Implementation}\label{gibbscluster}

For our motif clustering model, a Gibbs sampler could intuitively be
based on $p(\bTheta_i | \bTheta_{-i})$.  However, since our
$\bTheta_i$'s are actually unknown, a more efficient clustering
procedure involves drawing values of the clustering indicators
directly, as in \cite{Maceachern1994}, without dealing with drawing a frequency matrix $\bTheta_i$
for each motif $i$ at each iteration.  We denote our clustering indicators $\bz = (z_1, \ldots z_n)$, which simply defines a partition of $\{ 1,\ldots, n\}$.  Algorithmically, we let $z_i = c$ if
$\bTheta_i$ belongs to the $c$-th cluster or $z_i = 0$ if $\bTheta_i$ is drawn from
the prior measure $\gamma$,  and hence forms a new cluster.
Our ``collapsed" Gibbs sampler \cite[]{Liu1994} iteratively samples from 
$p(z_i | \bz_{-i},\bY)$
where we again use the notation $\bz_{-i}$ or $\bTheta_{-i}$ to mean all the
$\bz$ or $\bTheta$ parameters except the $i$-th one.  As mentioned earlier, no matter what the current $\bTheta_{-i}$ is, as long as they correspond to the same indicator vector $\bz_{-i}$, we have the same (almost surely) conditional prior distribution $p(z_i | \bTheta_{-i})$. Thus, we can write that
\begin{eqnarray}
p(z_i | \bTheta_{-i}) \equiv p (z_i = 0 | \bz_{-i})  =  \frac{b}{b+n-1} \hspace{2cm} p (z_i = c | \bz_{-i})  =  \frac{n_c}{b+n-1} \label{priorprobs}
\end{eqnarray}
where $n_c$ is the size of cluster $c$ (ie. the number of $z$'s in $\bz_{-i}$ which are equal to $c$) and $b$ is the weight parameter for forming a new cluster.  
Thus, under this model the prior probability for joining a particular cluster increases as the number of observations 
in that cluster increases, implying that the Dirichlet process prior favors unequal allocations of observations.
With observations $\bY$, we have the posterior conditional distribution 
\begin{eqnarray*}
p(z_i | \bz_{-i},\bY) \propto p(\bz, \bY) \propto p(\bY \mid \bz) p(z_i \mid \bz_{-i}).
\end{eqnarray*}
It is evident from our model that $p(\bY \mid \bz)$ is straightforward to derive and can be written as the product of the normalizing constants of $C$ product multinomial distributions, where $C$ is the total number of distinctive $z$Õs.  After some simplifications, we have the probability that observation $Y_i$ forms a new cluster is
\begin{eqnarray}
p(z_i = 0 | \bz_{-i}, \bY) & \propto & \frac{b}{b+n-1} \prod\limits_{j=1}^w 
\frac{\prod_k \Gamma(Y_{ijk} + \alpha)}{\Gamma(\sum_k Y_{ijk} + 4\alpha)} 
\frac{\Gamma(4\alpha)}{\Gamma(\alpha)^4} ,\label{newclusprob}
\end{eqnarray}
where $b$ is the weighting factor as defined at the end of Section~\ref{hierarchicalmodel}.  For the case where $z_i = c \neq 0$ (i.e., joining an existing cluster that already has 
a count matrix $\tilde{\bY}_c$), we have 
\begin{eqnarray}
\hspace{-1.25cm} p(z_i = c \mid \bz_{-i}, \bY) & \propto & \frac{n_c}{b+n-1}  \prod\limits_{j=1}^w 
\frac{\prod_k \Gamma(Y_{ijk} + \tilde{Y}_{cjk} + \alpha)}{\Gamma(\sum_k Y_{ijk} + \tilde{Y}_{ljk} + 4\alpha)} 
\frac{\Gamma(\sum_k \tilde{Y}_{cjk} + 4\alpha)}{\prod_k \Gamma(\tilde{Y}_{cjk} + \alpha)} \label{curclusprob}
\end{eqnarray}
A complete iteration of our Gibbs sampling algorithm results in a complete sample {\bf z} 
of our clustering indicators, which also represents a complete {\it partition} of our motif 
matrices. 

\subsection{Motif Alignment}\label{motifalignment}

An additional component of our model addresses the fact that  
we do not necessarily know which {\it core} $\bY_i$ of width $w$ to use 
within the raw alignment matrix of width $n_i > w$ for motif $i$.  We use the notation $a_i = j$ to mean that we are using 
the columns $j,j+1,j+w-1$ of our raw motif matrix $\bN_i$ as our core $\bY_i$.  For example, if 
our clustering algorithm is based on a fixed width of $w=6$ and our $i$-th raw motif matrix 
$\bN_i$ has 8 positions, then we have three possible choices for our core motif: $a_i = 1$ ($\bY_i$ = columns 
1 to 6 of $\bN_i$), $a_i = 2$ ($\bY_i$ = columns 2 to 7 of $\bN_i$), or $a_i = 3$ ($\bY_i$ = columns 3 to 8 of $\bN_i$).
We thus need an additional step where, for each raw data matrix, the  
best location of the central motif 
$a_i$ is drawn conditional the other motifs $\bY_{-i}$ and clustering indicators $\bz_{-i}$ for the other motifs.   Let $\bY^{a_i}_i$ denote the core $\bY_i$ that corresponds to the choice of a particular $a_i$, then the posterior probability $p(a_i) =  p(a_i|\bz_{-i},\bY_{-i})$ of $a_i$ is 
{\footnotesize 
\begin{eqnarray}
\hspace{-0.5cm} p(a_i) & \propto & \int p(\bY^{a_i}_i | \bTheta_i) \; 
p(\bTheta_i | \bz_{-i}, \bY_{-i}) \; d \bTheta_i \nonumber \\
& = & \int \bTheta_i^{\bY^{a_i}_i} \; \frac{b}{b+n-1} \; \left[\frac{\Gamma(4\alpha)}{\Gamma{\alpha}^4} 
\right]^w \bTheta_i^\alpha \; d \bTheta_i  + \sum\limits_{c=1}^{C}  \int \bTheta_i^{\bY^{a_i}_i} \; \frac{n_c}{b+n-1} \; 
\bdelta_{(\bTheta_i = \tilde{\bTheta}_c)} \; \frac{p(\tilde{\bY}_c,\tilde{\bTheta}_c | \bz_{-i})}{p(\tilde{\bY}_c | \bz_{-i})} \; d \bTheta_i \nonumber \\
& =  & \frac{b}{b+n-1} \left[\frac{\Gamma(4\alpha)}{\Gamma{\alpha}^4} 
\right]^w \frac{\prod_k \Gamma (Y^{a_i}_{ijk}+\alpha)} 
{\Gamma (\sum_k Y^{a_i}_{ijk} + 4\alpha)} 
+ \sum\limits_{c=1}^{C} \frac{n_c}{b+n-1}  \prod\limits_{j=1}^w 
\frac{\prod_k \Gamma(Y^{a_i}_{ijk} + \tilde{Y}_{cjk} + \alpha)}{\Gamma(\sum_k Y^{a_i}_{ijk} + \tilde{Y}_{cjk} + 4\alpha)} 
\frac{\Gamma(\sum_k \tilde{Y}_{cjk} + 4\alpha)}{\prod_k \Gamma(\tilde{Y}_{cjk} + \alpha)}
\end{eqnarray}
}
This alignment procedure is performed every tenth iteration of the collapsed Gibbs sampler described in the previous 
section.  

\subsection{Allowing Cluster Width to Vary}\label{varyingwidth}

We now extend our model to allow the core motif widths within each cluster, $\bw = (w_1, \ldots, w_C)$, where $C$ is the current number of clusters, to be unknown variables.   Each cluster width $w_c$ is modeled as being independent with prior distribution 
$w_c \sim {\rm Poisson} (\lambda)$
where $\lambda$ is the expected {\it a priori} width of the motif in each cluster, which we assume is fixed and specified.   
We let $B_{ik}$ be the ``background" counts of nucleotide $k$ over all columns of the raw matrix $\bN_i$ that are not included in the core matrix $\bY_i$, which now must be taken into account by our model since each motif width is allowed to vary.   These background columns represent the edges of the raw motif matrices that are not well conserved. We assume that the background counts,  
$\bB_i = (B_{ia}, \ldots B_{it})$, are a multinomial realization from an underlying vector of background nucleotide frequencies $\btheta_{0}= (\theta_{0a}, \ldots, \theta_{0t})$.
With these added distributions, the posterior probability for the core width of a particular cluster $c$, conditional on the current members of that cluster $\bz_c$ and their core alignments $\ba_c$, is 
\begin{eqnarray}
p(w_c \mid \bz_c, \ba_c) \propto 
\prod\limits_{j=1}^{w_c} \frac{\prod_k \Gamma(\tilde{Y}_{cjk} + \alpha)}{\Gamma(\sum_k  \tilde{Y}_{cjk} + 4\alpha)} \cdot
\frac{\Gamma(4\alpha)}{\Gamma(\alpha)^4} \cdot \prod_k \theta_{0k}^{B_{ck}} \cdot \frac{\lambda^{w_c} e^{-\lambda}}{\Gamma (w_c)} 
\label{widthposterior}
\end{eqnarray}
where $\tilde{Y}_{cjk}$ and $B_{ck}$ are, respectively, the core and background nucleotide counts in cluster $c$.   We implement this added component of our model by an additional step in our Gibbs sampling algorithm that, for each cluster $c$, samples a new value of 
$w_c$ from the conditional posterior distribution (\ref{widthposterior}).    Some methods, such as \cite{SchSumZha2005}, use a fixed number of columns to calculate the similarity between matrices and these core widths are then held fixed during their clustering procedure.  In contrast, our Gibbs sampling implementation allows the core width of each motif to change depending on the current set of clusters.  This means that each core width in our model is being estimated using additional information from other motifs in the dataset, rather than independently estimating each core width using only information from each matrix separately.  Other methods, such as \cite{Kielbasa2005}, cluster motifs on the basis of a pairwise distance and also estimate core width based on a pairwise comparison.   This is a slight improvement compared to treating each matrix individually, but still does not use as much information as our model, which makes width decisions on the level of an entire cluster, not just for pairs of motifs.  

\subsection{Alternative Clustering Priors}\label{comparingpriors}

As mentioned in Section~\ref{gibbscluster}, the Dirichlet process
prior favors unequal allocation of observations, meaning that each new
observation has a greater prior probability of being placed in a
cluster that already has many observations.  
An alternative is a uniform clustering prior which favors equal allocations 
of observations i.e., the prior probability that 
a new observation is placed in any one of the existing clusters is uniform.
If we already have $n$
observations divided into clusters $c=1,\ldots,C$ with $n_1,\ldots,n_c$
members, then the two prior distributions are

\begin{center}
\begin{tabular}{lll}
{\underline{Dirichlet process prior}} & \hspace{2cm} & {\underline{uniform clustering prior}} \\
$P(z_{n+1} = c | \bz, {\rm DP}) = \frac{n_c}{b+n}$ & & $P(z_{n+1} = c | \bz, {\rm unif}) = \frac{1}{b+C} $ \\ 
$P(z_{n+1} = {\rm new} | \bz, {\rm DP}) = \frac{b}{b+n} $ &  & $P(z_{n+1} = {\rm new} | \bz, {\rm unif}) = \frac{b}{b+C} $\\
\end{tabular}
\end{center}

where $b$ is the weight given to forming a new cluster.  In fact, we can consider both the Dirichlet process and uniform
clustering specifications as particular cases of a more general
clustering prior distribution, where
\begin{eqnarray}
P(z_{n+1} = {\rm new} | \bz)  \propto b \hspace{2cm} P(z_{n+1} = c | \bz)  \propto  f(n_c) \quad\quad c = 1, \ldots, C \label{generalpriorprobs}
\end{eqnarray}
This general clustering model reduces to the Dirichlet
process when $f(n_c) = n_c$ and the uniform clustering prior when
$f(n_c) = 1$, but more general functions may be desirable in
particular situations.  The prior density of a partition $\bz$ with $C$ clusters under either our Dirichlet
process or uniform clustering model can be calculated recursively. With the 
Dirichlet process prior, we have prior density 
\begin{eqnarray}
p(\bz | {\rm DP}) = \frac{b^C \cdot \prod\limits_{c=1}^{C} (n_c - 1)!}{\prod\limits_{i=1}^n (b + i - 1)} \label{DPpriordensity}
\end{eqnarray}
With the uniform clustering prior, we have prior density
\begin{eqnarray}
p(\bz | {\rm unif}) = \frac{b^{C-1} \cdot (b+C)}{\prod\limits_{c=1}^C (b + c)^{n_c}} \label{Unifpriordensity}
\end{eqnarray}

It is worth noting that the Dirichlet process prior density (\ref{DPpriordensity})  is not affected by the ordering of the clusters, whereas the uniform clustering prior density (\ref{Unifpriordensity}) is.  In other words, different partitions
with the same cluster sizes are exchangeable under the Dirichlet
process model, but we will  get different values of the uniform clustering density 
for different, but exchangeable orderings of unequally-sized
clusters. As suggested by \cite{GreRic2001}, to ensure exchangeability of our
uniform clustering model, we need to make our prior density $p(\bz | {\rm
Unif})$ a function of a ``signature" of the partition that is
identical for exchangeable partitions.  For example, if we let $p^{\ast}(\bz |
{\rm Unif}) = k \cdot p(\bz^{\prime} | {\rm Unif})$ where
$\bz^{\prime}$ is $\bz$ with the $z_i$'s arranged in order from the
largest cluster to the smallest, then the calculation of
(\ref{Unifpriordensity}) for $\bz^{\prime}$ will be the same for all
exchangeable values of $\bz$.  All of these complications are avoided
in the Dirichlet process model which automatically gives the same
prior density value for exchangeable partitions.

We compared the behavior of these two clustering prior
specifications for a simple simulation study, where 1000 complete
partitions $\bz = (z_1,\ldots,z_n)$ with $n = 106$ and $b = 1$ were
generated under both prior distributions.  This particular sample size $n$ was chosen to be comparable with the real data in our first application in Section~\ref{jasparapp} below.  Figure~\ref{priorclusinfo} displays the
distributions of both the number of clusters as well as the size of
the multiple-member ($n_c > 1$) clusters over all of our simulated
partitions.  As expected, the number of clusters (with multiple members) is larger under the uniform 
prior and the size of some clusters from the Dirichlet process are larger than any 
generated from the uniform prior specification.  

\begin{figure}[ht]
\caption{Comparison of clustering statistics between DP and uniform 
clustering priors}\label{priorclusinfo}
\begin{center}
\rotatebox{270}{\includegraphics[width=3in,height=5in]{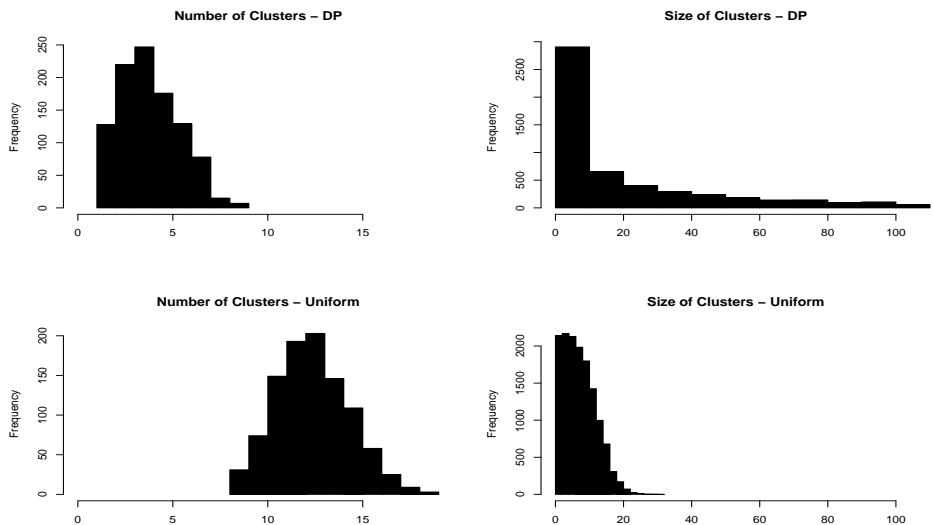}}
\end{center}
\end{figure}

We now discuss several applications of our Bayesian motif clustering model to collections of transcription factor matrices.  In Section~\ref{jasparapp}, we apply our method to JASPAR \cite[]{jaspar}, a small but heavily curated database of transcription factor matrices.   In Section~\ref{combined}, we apply our method to a collection consisting of all matrices from JASPAR and the TRANSFAC database \cite[]{transfac}, which is a larger but less curated database of transcription factor matrices.   Finally, in Section~\ref{crossclus}, we discuss a cross-species application where discovered motifs are used to infer co-regulated genes in bacteria.  

\section{Application to a single database: JASPAR} \label{jasparapp}

We use the JASPAR \cite[]{jaspar} database as an example for illustrating different strategies for visualizing and analyzing results from our clustering model. This database contains 111 nucleotide-count matrices which differ substantially in appearance,
number of counts, and motif width.  
Between different motifs, the number of binding sites used to construct the motif matrix (the total number of multinomial counts) varies from 6 to 389, with an average of around 35 counts per matrix.  The range of matrix widths was from 4 to 30 bp, though the matrices were generally short, with an average width of approximately 11 bps. 
From this database, we also have species information for almost all motifs, as well as a 
classification into a
particular ``protein family'' based on the common physical structure
of each motif's DNA-binding domains.  
For example, one family of
transcription factors is the helix-loop-helix family, which has two
DNA-binding helix domains that bind directly to the DNA strand and are
joined by a loop domain.   We will use this extra protein family and species information when we examine the results produced by our clustering model.

For our
Dirichlet process prior, we chose prior parameters $\alpha$ and prior weight $b$ to both be equal to 1.  We chose a small prior weight $b$ for simplicity, since we had no prior knowledge on the number of clusters to expect, and our clustering results did not change substantially when using larger values of $b$.  We also assumed our background frequencies were uniform across nucleotides, ie. $\theta_{0k} = 0.25$ for $k = 1, \ldots, 4$.  Although we allowed the motif widths to vary, we restricted the core width of each motif to be at least 6 bps in order to reflect biological reality, with an {\it a priori} expected core width 
$\lambda$ of 8 bps.  This restriction reduced our dataset from 111 matrices down to 106 matrices.  As described in Sections~\ref{gibbscluster}-\ref{varyingwidth}, our Bayesian hierarchical
clustering model was implemented using a Gibbs sampling algorithm.  Details of our evaluation of convergence are given in the supplemental materials.  

\subsection{Tree of Pairwise Clustering Probabilities}\label{clustrees}

\begin{figure}[ht]
\caption{Tree of Pairwise Clustering Probabilities.  The motifs are labeled by both the motif name as well as the protein
family and species associated with that motif. Motifs which did not cluster
with any other motifs are not shown.  The length of each tree
``branch'' (vertical line shared by a group of motifs) is proportional to the average distance $d_{ij}$ between that group of motifs and the next closest group of motifs. The arbitrary ordering of the branches is 
chosen by the hierarchical algorithm.}  \label{variablewidthtree}
\vspace{-1cm}
\begin{center}
\hspace{-1.5cm}
\includegraphics[width=7.5in,height=7.5in]{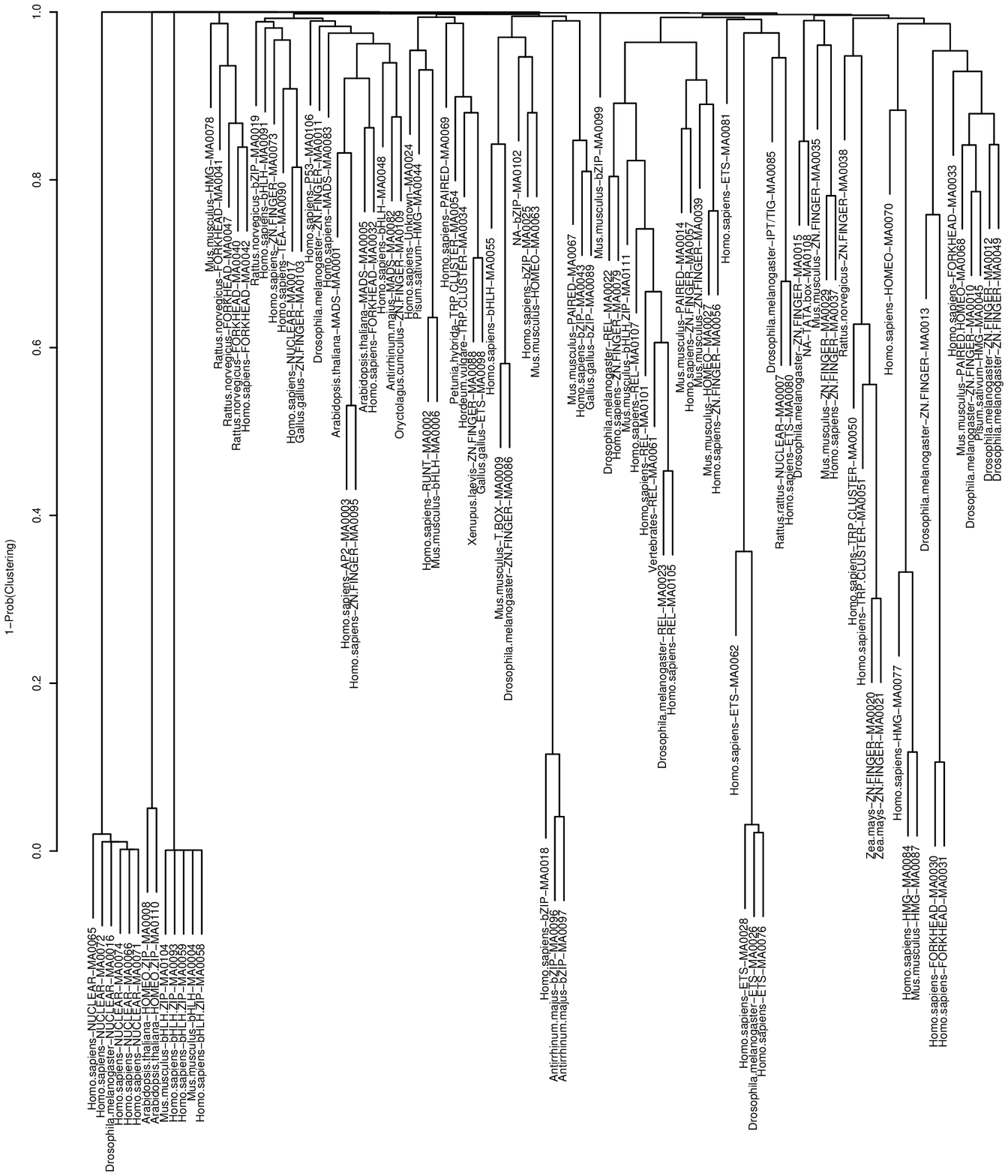}
\end{center}
\end{figure}

An intuitive means for examining our overall clustering results is the
posterior probability $p_{ij}$ that a particular pair of motifs $i$
and $j$ are in the same cluster. The value of $p_{ij}$ for any two
motifs $i$ and $j$ can be estimated by the proportion of iterations
that have motif $i$ and $j$ in the same cluster, which is the 
Monte Carlo estimate of the posterior mean of the indicator variable
for motif $i$ and $j$ being in the same cluster.  The full list of pairwise clustering probabilities is given in the supplemental materials.  Based on these pairwise clustering probabilities $p_{ij}$, a pairwise
distance measure can be calculated between each pair of motifs in the
dataset, $d_{ij} = 1-p_{ij}$.  
 \cite{MedSiv2002} use the same distance measure in their gene expression application.  Our distance matrix was then converted into a tree diagram by an average-linkage hierarchical algorithm, which is 
shown in Figure~\ref{variablewidthtree}. 


This clustering tree shows several strong relationships, such as the group of NUCLEAR motifs and 
the group of bHLH motifs on the lefthand side of Figure~\ref{variablewidthtree}.  
Many weaker relationships are also present, implying that many motifs have a low but non-zero probability of
being grouped together. Although many of the stronger clusters of motifs belong to the same protein family and same species,
there are several interesting exceptions.   The NUCLEAR and bHLH groups mentioned above contain mostly motifs from {\it Homo sapiens}, but the NUCLEAR group also contains a motif from fruit fly {\it Drosophila melanogaster} and the bHLH group includes motifs from mouse ({\it Mus musculus}).   
                                                                                         
\subsection{Best Clustering Partition and Cluster Strength}\label{bestpartition}

Although Figure~\ref{variablewidthtree} allows us to examine the clustering structure of the
entire dataset, the tree is not ideal for deducing the
``best partition" or best set of clusters in the dataset, since these pairwise probabilities are calculated across many different
partitions, similar to the problem for  hierarchical tree clustering mentioned in Section~\ref{clusintro}. 
One could ``cut the tree" at any number of different
threshold distances and thereby produce any number of possible
partitions, but a less arbitrary alternative is to estimate of the posterior mode of our clusters from our MCMC simulation. 
We estimate this
posterior mode by calculating the posterior probability of the partition
$\bz$ at the end of each iteration of our sampler, and retaining the
partition $\hat{\bz}$ with the highest posterior probability as our best
estimate of the mode.
The probability $p(\bz | \bY, \bB)$ of $\bz$ is calculated as the
product of the likelihood value of our cluster matrices $\tilde{\bY}_c$ (the sum of all $w_c \times 4$
count matrices in cluster $c$) and total background counts $B_k$, 
\begin{eqnarray*}
p(\bY,\bB|\bz) \quad \propto \quad \prod\limits_{c=1}^{C} \int p(\tilde{\bY}_c |\bTheta_c)
p(\bTheta_c|\bz) p(\bB | \btheta_0) d \bTheta_c 
\quad \propto \quad \prod\limits_{c=1}^{C} \prod\limits_{j=1}^{w_c} \frac{\prod_k \Gamma(\tilde{Y}_{cjk} + \alpha)}{\Gamma (\sum_k \tilde{Y}_{cjk} + 4\alpha)} \prod\limits_k \theta_{0k}^{B_k} 
\end{eqnarray*} 
and the prior densities of our clusters $p(\bz)$ and variable motif widths $p(\bw)$.  
Although this best partition can reduce our dataset down to a list of
interesting clusters, we have lost information about the variability
of these clusters by focusing on a point estimate. 
In order to retain some measure of variability, we incorporate cluster-level and observation-level clustering
characteristics within our ``best clusters". We
can measure the strength of each cluster by calculating the logarithm of the {\it Bayes
factor} \cite[]{KasRaf1995} for the current cluster $c$, with members
$\bz = (z_1,z_2,\ldots,z_{n_c})$, versus each member of the cluster
forming its own cluster,
\begin{eqnarray*}
{\rm Strength}({\rm Cluster} \; c) = \log \left[
\frac{P(\bY | \bz \; {\rm all} \; {\rm same} )}{P(\bY | \bz \; {\rm
all} \; {\rm different} )} \times \frac{P(\bz \; {\rm all} \; {\rm
same} )}{P(\bz \; {\rm all} \; {\rm different} )} \right]
\end{eqnarray*} 
For a cluster of motifs ($\bY_1,\ldots,\bY_m$) and clustering
indicators $\bz = (z_1,\ldots,z_m)$,
\begin{eqnarray*}
{\rm Strength} = \log \left[ \frac{\int \tilde{\bTheta}^{\tilde{\bY} + \alpha
- 1} d \tilde{\bTheta}}{ \prod\limits_{i=1}^{m}
\int \bTheta_i^{\bY_i + \alpha - 1} d \bTheta_i} \times \frac{(m-1)!}{b^{m-1}}  
\right] 
= \log \left[ \frac{\prod\limits_{j=1}^{w_c} \frac{\prod_k \Gamma (\tilde{Y}_{jk} +
\alpha)}{\sum_k \Gamma (\tilde{Y}_{jk} + 4 \alpha)}}{\prod\limits_{i=1}^m
\prod\limits_{j=1}^{w_c} \frac{\prod_k \Gamma (Y_{ijk} + \alpha)}{\sum_k \Gamma(Y_{ijk} + 4 \alpha)}}
\times \frac{(m-1)!}{b^{m-1}} \right]
\end{eqnarray*}
where $\tilde{\bY}$ and $\tilde{\bTheta}$ again denote the count and
frequency matrices for the entire cluster together.  
The clusters within our best partition can then be ranked by this
measure of cluster strength, giving us an extra measure of
confidence/uncertainty about inference based upon a specific cluster.
We can also measure clustering strength at the level of individual
motifs within our best partition by calculating, for each motif, the
posterior probability that it should belong to that cluster, as
opposed to any of the other existing clusters or being its own
cluster.  For each motif $i$, this posterior probability $p(z_i |
\bz_{-i}, \bY)$ is the same calculation that is performed during each
iteration of our Gibbs sampling algorithm, but in this case we are
conditioning on the best partition i.e., $p(\hat{z}_i | \hat{\bz}_{-i},
\bY)$.

For our dataset, the partition $\hat{\bz}$ with the highest posterior value consisted of 
26 multiple-member clusters containing 69 out of 106 total motifs.  The strongest 6 of these 26 clusters (using our cluster strength measure) are listed in Table~\ref{bestpartitiontable}, along with their cluster size, width and consensus sequence.  The consensus sequence is a 
representation of the total count matrix for the cluster, giving the nucleotide with the highest 
count in each position. A nucleotide is only capitalized if its nucleotide frequency is greater 
than 0.75 in that position.  It is clear from the
table that this measure of
cluster strength is quite dependent upon the size of the cluster:
larger clusters tend to have a higher value of cluster strength.

\begin{table}[ht]
\caption{Strongest Clusters from the Best Partition of JASPAR matrices}\label{bestpartitiontable}
\begin{center}
\begin{tabular}{|ccccccll|}
\hline
Clus & Size & Strength & Width & Consensus & Protein Families & Species & Motifs \\ 
\hline
1 & 6 & 186.5 & 6 & {\tt aGGTCA} & NUCLEAR & {\it H. sapiens} &  MA0065 MA0066  \\
  &  &   &   &   &   &  & MA0071 MA0072 \\
  &  &   &   &   &   &  &  MA0074 \\
  &  &    &   &   &   &   {\it D. melanogaster}  & MA0016 \\
\hline
2 & 5 & 179.5 & 6 & {\tt CACGTG} & bHLH-ZIP & {\it H. sapiens} & MA0058 MA0059   \\
  &  &   &   &   &   &  &  MA0093 \\
&  &   &   &   &   & {\it M. musculus} & MA0104 \\
 &  &   &   &   &  bHLH & {\it M. musculus} & MA0004 \\
\hline
3& 3 & 82.2 & 7 & {\tt aTGACGT} & bZIP & {\it A. majus} & MA0096 MA0097 \\
 &   &   &   &   &    & {\it H. sapiens} & MA0018 \\
\hline
4 & 3 & 72.7 & 7 & {\tt CCGGAAg} & ETS & {\it H. sapiens} & MA0028 MA0076 \\
 &   &   &   &   &    & {\it D. melanogaster} & MA0026 \\
\hline
5 & 3 & 72.5 & 8 & {\tt GTAAACAa} & FORKHEAD & {\it H. sapiens} & MA0030 MA0031 \\
 &   &   &   &   &  ZN-FINGER  & {\it D. melanogaster} & MA0013 \\
\hline
6 & 2 & 40.0 & 8 & {\tt cAATtATT} & HOMEO-ZIP & {\it A. thaliana} & MA0008 \\
 &   &   &   &   &  & {\it A. thaliana} & MA0110 \\
\hline
\end{tabular}
\end{center}
\end{table}

Most of our strongest clusters contain motifs from within a single TF protein family, though there are exceptions, such as the fifth cluster that contains both ZN-FINGER and FORKHEAD motifs.  Clearly, this cluster would not have been detected if motifs were only grouped together based on TF family.  
It is also interesting to note that most of the larger clusters contain similar motifs from different species, with motifs from human ({\it H. sapiens}) being grouped with motifs from mouse ({\it M.musculus}), fruit fly ({\it D. melanogaster}) and snapdragon flowers ({\it Antirrhinum majus}).  
All of the motifs in Table~\ref{bestpartitiontable} have individual clustering probabilities close to 1, but many of the weaker clusters in the best partition contain motifs which have individual probabilities substantially less than one.  Our entire best partition of JASPAR matrices is given in the supplemental materials.  

\subsection{Core Width Variability in JASPAR}\label{clusteringwidth}

\begin{figure}[ht]
\caption{95\% posterior intervals of each motif width in JASPAR matrices} \label{widthintervals}
\begin{center}
\rotatebox{270}{
\includegraphics[width=3.5in,height=6in]{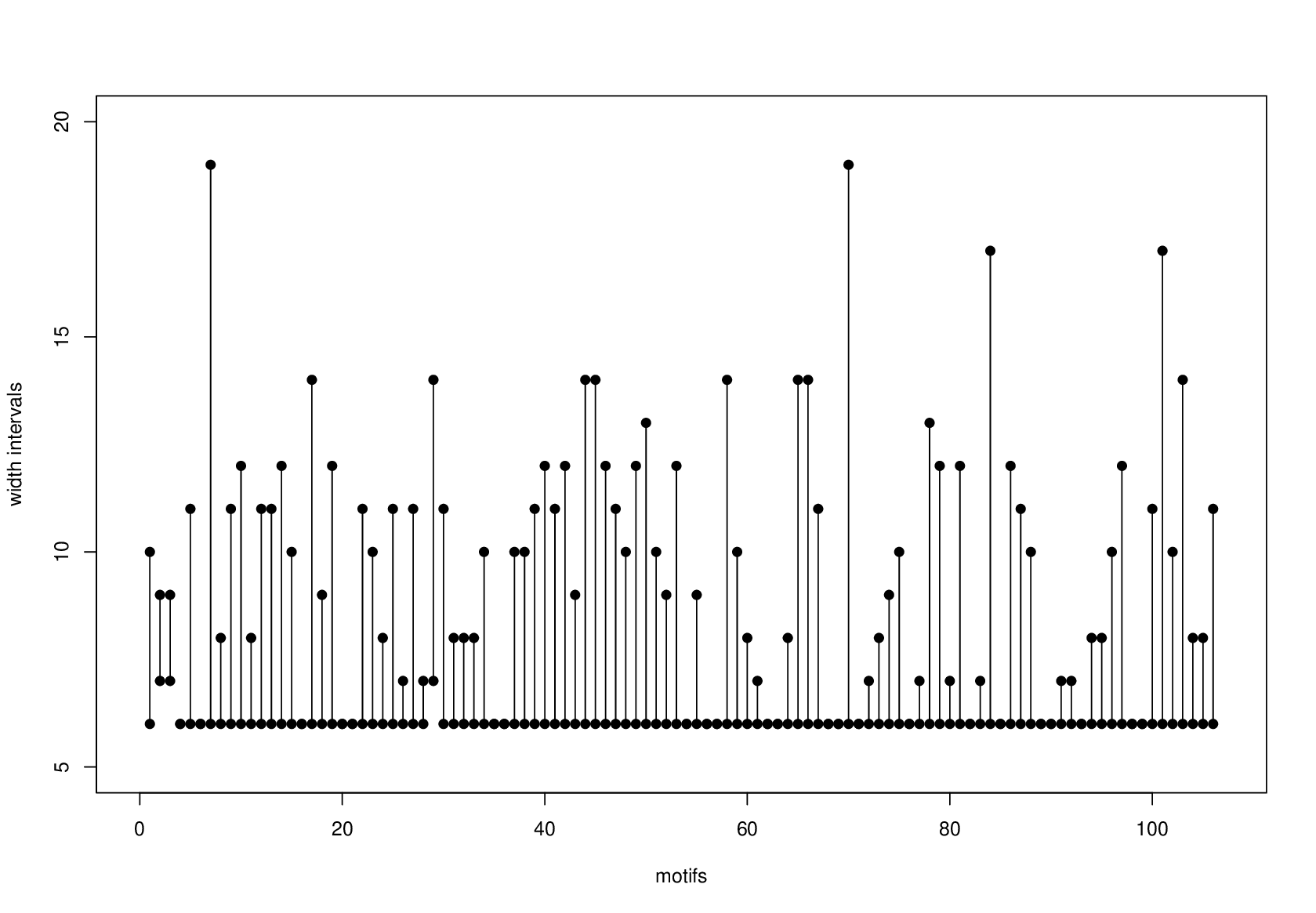}
}
\end{center}

\vspace{-0.75cm}

\end{figure}

A key component of our clustering procedure is that the width of the core motif within each raw matrix is not considered to be fixed and known, but is instead allowed to vary by cluster in our model. We examine the variability of these core widths in Figure~\ref{widthintervals},  which is a plot of  95\% posterior intervals for the core width of each motif in our JASPAR dataset.  

The 95\% posterior intervals are quite different between motifs, with some motifs having a wide interval while other motifs have an interval consisting only of the minimum motif width of 6 bps, which in many cases is because the raw data matrix is only 6 bps wide.  In several motifs, the 95\% posterior interval for the motif width does not even include the {\it a priori} expected motif width of 8 bps.  The wide intervals for several motifs suggests that considering core widths as fixed and known, as in \cite{matinspector} and \cite{SchSumZha2005}, would result in a substantial loss of information

\subsection{Effect of Prior Specification on JASPAR results}

\begin{figure}[ht]
\caption{Clustering statistics between uniform and DP models} 
\label{chipdpvsuniformstats}
\vspace{-0.75cm}
\begin{center}
\rotatebox{270}{
\includegraphics[width=4in,height=6in]{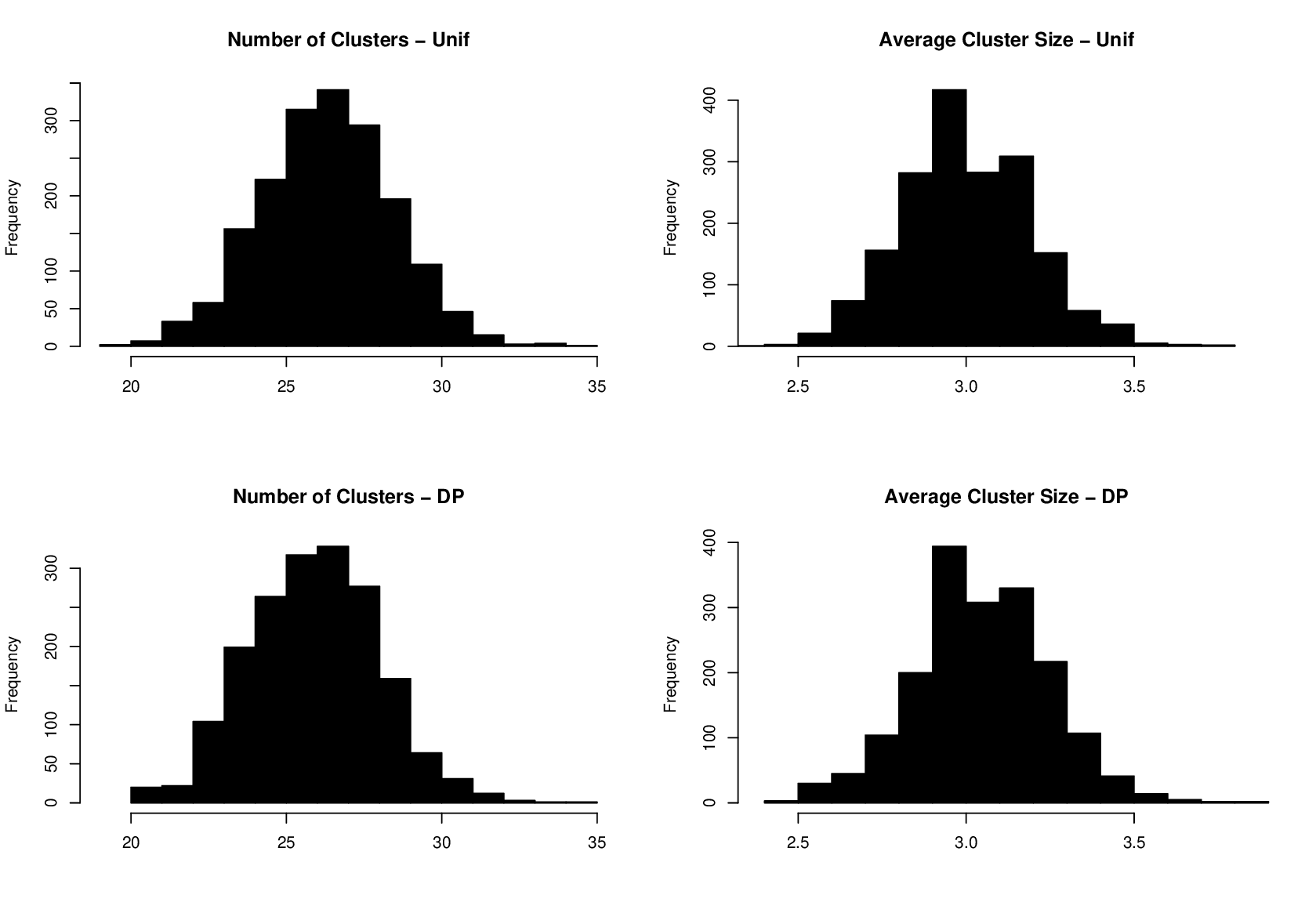}
}
\end{center}

\end{figure}

In Section~\ref{comparingpriors}, we discussed the {\it a priori} differences between the Dirichlet process prior
compared to the uniform clustering prior.  We now investigate whether these differences are also apparent in 
the posterior JASPAR results.  The distribution (over all partitions produced by the Gibbs sampler) 
of the number of multiple-member clusters and the average size of our JASPAR
clusters is given in Figure~\ref{chipdpvsuniformstats}.   

The Dirichlet process prior and uniform prior give dramatically different clustering results 
based upon prior simulation alone (Figure~\ref{priorclusinfo}) , but show very slight differences in the posterior clustering 
results (Figure~\ref{chipdpvsuniformstats}).   Only minor differences were observed between our Dirichlet process and uniform clustering models in terms of 
the clustering trees (Section~\ref{clustrees}) and best partitions (Section~\ref{bestpartition}), which indicates that our 
choice of a Dirichlet process prior distribution was not very influential on our posterior clustering results, at least in comparison to a uniform clustering prior alternative.  However, other datasets may show a larger influence of the 
prior specification on the posterior clustering results, and our software allows the user to specify the use of either prior distribution.  \cite{GreRic2001} demonstrate with several 
datasets that the unequal allocations favored by the Dirichlet process priors can persist in the 
posterior distribution.  

\section{Application to Combined Databases: JASPAR and TRANSFAC} \label{combined}

Our motivations for applying a clustering model to established databases, such as JASPAR or TRANSFAC, are to eliminate ``redundant" matrices within these databases and to understand the similarities as well as differences among the different transcription factors.  If two (or more) motifs have nearly identical core matrices, we can consider them as redundant since they do not contain unique information in terms of appearance.  The presence of redundant matrices complicates the use of these databases for further motif discovery, since searches involving several nearly identical matrices will lead to an excessive number of predicted sites \cite[]{Kielbasa2005}, and it also affects how we evaluate computational results.  As mentioned in Sections~\ref{motifalignment}-\ref{varyingwidth}, an additional complication in these databases is that the motif core has unknown width and location within each matrix, which is also addressed by our clustering model.    
In Section~\ref{jasparapp}, we already observed substantial redundancy within the JASPAR database alone.  Now, we want to combine our JASPAR matrices together with the larger TRANSFAC database \cite[]{transfac}, which contains 714 transcription factor matrices.   In this combined collection of 825 matrices, we expect substantial redundancy which we will address with our clustering procedure.  Similar to our JASPAR-only application, we chose prior parameters $\alpha$ and $b$ to both be equal to 1,  expected core width $\lambda$ = 8 bps, and background frequencies $\theta_{0k} = 0.25$ for $k = 1, \ldots, 4$. We again restricted the core width of each motif to be at least 6 bps, which reduced our dataset from 825 matrices down to 817 matrices.  Evaluation of Gibbs sampling convergence for this combined dataset are given in the supplemental materials.    

The best partition produced by our Gibbs sampling implementation consisted of 165 multiple-member clusters containing 746 out of the 817 matrices in our combined collection.  The clusters in our best partition are given in the supplemental materials, along with a full list of the pairwise clustering probabilities. We also examined this combined dataset using our model with the alternative uniform clustering prior, and  observed very little difference in the posterior clustering results (details given in supplemental materials).  Our expectations of high levels of redundancy in these databases seems confirmed by the fact that approximately 90\% of the matrices in our combined collection are partitioned into multiple-member clusters.  The strongest 10 of these 165 clusters (using our cluster strength measure) are listed in Table~\ref{bestpartitiontable2}, along with their cluster size, width and consensus sequence (as defined in our previous application).   For brevity, the individual members of each cluster are not given in Table~\ref{bestpartitiontable2}. Instead, the total number of members from either JASPAR or TRANSFAC in each cluster is given.

\begin{table}[ht]
\caption{Strongest Clusters from the Best Partition of TRANSFAC and JASPAR matrices}\label{bestpartitiontable2}
\begin{center}
\begin{tabular}{|ccccccc|}
\hline
Cluster & Size & Strength & Width & Consensus & \multicolumn{2}{c|}{Number of Members} \\
   &     &  &   &   &  JASPAR & TRANSFAC \\ 
\hline
1 & 24 & 1015.3 & 6 & {\tt CACGTG} &  1 & 23   \\
2 & 25 & 974.9 & 6 & {\tt TGACGT} &  2 & 23   \\
3 & 15 & 697.5 & 8 & {\tt TTTcGCGC} &  1 & 14   \\
4 & 18 & 550.1 & 6 & {\tt aGATAa} &  1 & 17   \\
5 & 16 & 548.4 & 6 & {\tt CACGTG} &  4 & 12   \\
6 & 17 & 486.5 & 6 & {\tt cGGAAg} & 3 & 14  \\
7 & 7 & 358.6 & 6 & {\tt TGTTCT} & 0 & 7 \\
8 & 11 & 338.7 & 7 & {\tt tcACGTG} & 0 & 11 \\
9 & 14 & 336.7 & 6 & {\tt TGAcCt} & 1 & 13 \\
10 & 11 & 319.3 & 6 & {\tt AAAGcg} & 4 & 7 \\
\hline
\end{tabular}
\end{center}
\end{table}

From Table~\ref{bestpartitiontable2}, we can see that there is substantial redundancy both between databases and within each database.  The redundancy seems to be more substantial within the TRANSFAC database, even taking into account the fact that the JASPAR matrices are a small component (13\%) of the combined collection.   This is not surprising, since there was a substantial amount of manual curation involved in the creation of the JASPAR database.   We also attempted additional merge moves to try and combine our strong clusters with similar appearance but these moves did not actually increase our posterior density.    Just as in our smaller JASPAR application (Section~\ref{jasparapp}), we again observed that a substantial number of motifs in our combined dataset also had substantial variability in their core widths and alignments.  87\% of the motifs in the combined dataset had 95\% posterior intervals for the core width which covered more than a single fixed value (details given in the supplemental materials).  It is also worth noting that a majority of the remaining 13\% of motifs had no variability simply due to the fact that the raw motif width was equal to the minimum core width.  Clearly, these results suggest that the assumption of fixed and known motif widths in \cite{matinspector} and \cite{SchSumZha2005} is a tenuous one, and results in the loss of substantial information.

Several existing motif-finding programs, such as MatInspector \cite[]{matinspector}, attempt to utilize non-redundant collections of matrices for scanning large sets of genomic sequences for transcription factor binding sites.  MatInspector uses a collection of matrix clusters created by manual curation, whereas we obviously prefer an automated method for generating clusters of matrices.  \cite{SchSumZha2005} present a clustering method for motif matrices which allows matrices to belong to multiple clusters, which means their method does not produce a {\it partition} of matrices which can be used to reduce redundancy.   
Our model implementation produces a best partition of clusters, alignments and core widths, all of which are needed to create a {\it super-matrix} for each cluster, which is the sum of the aligned core matrices in that cluster.  The set of super-matrices can then be used as inputs for a sequence scanning algorithm, such as \cite{Huangetal2004}.  
When looking for matches within genomic sequences, using these super-matrices reduces the number of comparisons that are needed compared to using each matrix in a database individually.   This reduction of redundant comparisons not only eases the computational burden of sequence scanning, but also helps somewhat with the usual problem of evaluating the statistical significance of good matches in this multiple comparison setting, as mentioned in \cite{Kielbasa2005}.

\section{Application to Cross-Species Conservation} \label{crossclus}

An additional application of our motif clustering procedure is presented in 
\cite{JenSheLiu2005},  where {\it phylogenetically}-discovered
motifs are clustered to infer co-regulated genes.   {\it Phylogenetic} motif discovery \cite[]{MccThomCarRyaLiuDerLaw2001} searches for conserved motifs in upstream sequences from different, but related, species under the assumption that transcription factor binding motifs are likely to be conserved by evolution.  If the motifs found 
upstream of several {\it Bacillus subtilis} genes are similar enough to be
clustered together, then it is possible that the same TF (recognizing
that common motif) is targeting each of the genes in that cluster.
Thus, by combining statistical techniques for both motif discovery and
motif clustering, one can infer clusters of potentially co-regulated genes.  \cite{QinMccThoMayLawLiu2003} inferred co-regulated genes by applying a clustering algorithm to previously discovered motifs in {\it E.coli}.  Instead of using a Dirichlet process prior,  \cite{QinMccThoMayLawLiu2003} used the same uniform clustering prior presented in Section~\ref{comparingpriors}.  However, their clustering techinique was simpler  than our model in that the core motif width was considered to be fixed and known, which we suggest in Sections~\ref{jasparapp} and~\ref{combined} is quite a restrictive assumption.  

In \cite{JenSheLiu2005}, we use phylogenetic motif discovery on genomic sequences from the bacteria {\it Bacillus subtilis} and six related bacterial species.  Our Bayesian hierarchical clustering model is then used to create clusters of highly-similar motifs within our collection of discovered motif matrices.  Since each motif contains sites discovered in close proximity to a particular {\it B.subtilis} gene, our motif clusters can be also interpreted as clusters of possibly co-regulated genes.   These gene clusters show substantial evidence of co-regulation when compared to several sources of external information about {\it B.subtilis} genes, such as microarray gene expression data and functional ontology.  In addition to using the best partition of clusters, we were also able to incorporate our clustering variability in this application.   Our cluster strength measure was again used to rank clusters, and several of the strongest clusters were also shown to be strong in terms of external evidence of co-regulation.  Our individual motif strengths $p(\hat{z}_i | \hat{\bz}_{-i}, \bY)$ were also used to filter out weakly-associated motifs from our clustering results.  

\section{Discussion}\label{discussion}

Although we have presented our clustering procedure in the context of 
several specific applications to motif matrices, 
these advantages of our Bayesian clustering model 
are not specific to this particular type of data.  Bayesian hierarchical clustering models 
based on a Dirichlet process prior distribution should be considered an attractive approach to many clustering problems.  
Our hierarchical framework lets us account for uncertainty in the
count matrices that represent each TF motif by assuming a product
multinomial distribution, whereas most clustering programs assume
the count matrices are fixed and known without error.
The second advantage is that our clustering strategy uses a model-based similarity 
measure rather than some ad hoc distance measure in order to compare motifs.  At each
iteration of the Gibbs sampling algorithm, the decision to cluster a
particular observation is determined by the conditional distribution
of $z_i$ given all other information ($\bz_{-i},\bY_i$).
Thus, our distance metric is exactly equal to the conditional
posterior distribution under our full Bayesian clustering model.
A third advantage of our clustering model allows not only the clusters
themselves to vary (in terms of which motifs are members of which
clusters) but also the number of clusters is allowed to vary.  This is
a key improvement over a clustering technique that requires the number of clusters 
to be fixed (such as K-means clustering) since, in this situation,  we have 
very little idea {\it a priori} about how many motifs we might expect would 
be similar to each other.  Standard hierarchical tree clustering is also less ideal in this situation, since an arbitrary threshold must be used to produce a set of clusters (eg. \cite{Kielbasa2005} and \cite{SchSumZha2005}). 
Another general advantage of our procedure is that our posterior sampling 
implementation gives us an idea of the variability of our
clustering results, whereas traditional clustering methods typically give only 
a point estimate.   We explore several summaries of this variability, including tree structures that summarize the pairwise clustering probability of our motifs, as well as measures of strength for entire clusters and individual motifs within clusters.  However, further research is needed into effective techniques for analyzing stochastic clustering results, since the usual procedure of averaging across iterations is not appropriate when both cluster sizes and individual memberships within clusters vary between iterations.  

We also presented a novel extension of our model that allows the motif width within each cluster to vary, and our results indicate that many motifs have substantial motif width variability.  Previous methods, such as \cite{QinMccThoMayLawLiu2003} may be ignoring important information by considering motif core widths to be fixed and known {\it a priori}.  Our model also addresses the alignment issue that, within each raw motif matrix, it is not obvious where the core motif is located.   Our model allows us to condition on the  
motif core in all other raw matrices within the current cluster when we calculate the most 
likely location of the motif core within a particular matrix.  In many cases, other matrices 
may show very similar compositions to the matrix in question (especially matrices within the same cluster), in which case the 
conditioning provides a substantial amount of information pertaining to the motif 
core location.    This extra information is ignored by methods (eg. \cite{SchSumZha2005}) which use fixed widths during their clustering procedure, and only partially captured by methods (eg. \cite{Kielbasa2005}) which use widths estimated by pairwise comparisons of matrices.

As mentioned in Section~\ref{combined}, our clustering results eliminate the redundancy within current matrix databases, which will benefit future motif discovery by reducing the number of redundant hits when scanning sequences for known transcription factors.  However, it may be possible to utilize our clustering model for motif discovery in a more sophisticated way.  A Bayesian framework for motif discovery
was presented in \cite{JenLiuZhoLiu2004} based on a motif
model where very little is known {\it a priori} about the
appearance of an unknown motif.  
However, once a set of motifs has
been discovered (and clustered), we should incorporate this
information directly into our motif discovery procedure.  
One proposal would be to
use the posterior predictive distribution from our motif clustering
model as the scoring function for motif discovery, which would
increase the ability of our motif-finding algorithms to detect a motif
that is similar to motifs that have already been discovered elsewhere.  In addition to using our best partition to aide the discovery of new transcription factor binding sites, this strategy would also allow us to utilize the clustering uncertainty and variability in motif core widths and alignments which is estimated by our model.  

\section*{Acknowledgment}

The authors thank Geetu Tuteja for help with the JASPAR and TRANSFAC databases. This research is partially supported by NSF DMS-0204674 and NSFC
10228102/A010201. 

\bibliography{references}

\end{document}